\documentclass[10pt]{article}
\usepackage{amssymb}
\usepackage{amsmath}
\numberwithin{equation}{section}

\newtheorem{theorem}{Theorem}

\newtheorem{lemma}{Lemma}

\newcommand{\ooo}{\overline}
\newcommand{\ep}{\varepsilon}

\newcommand{\ppp}{\partial}
\newcommand{\www}{\widetilde}
\newcommand{\OOO}{\Omega}

\numberwithin{equation}{section}

\newcommand{\ddda}{d_t^{\alpha}}

\newcommand{\LLLLL}{L^2(0,\ell)}

\allowdisplaybreaks

\title
{\bf Inverse coefficient problems for one-dimensional time-fractional
diffusion equations}

\author{$^1$ Oleg Imanuvilov, $^2$ Kazufumi Ito and  
$^{3,4}$ Masahiro Yamamoto}

\date{}
\begin{document}
\maketitle

\thanks{
$^1$ Department of Mathematics, Colorado State University, 101 Weber Building, 
Fort Collins CO 80523-1874, USA 
e-mail: {\tt oleg@math.colostate.edu}
$^2$ Department of Mathematics, North Carolina State University, 27695,
USA e-mail:{\tt kito@ncsu.edu}
$^3$ Graduate School of Mathematical Sciences, The University
of Tokyo, Komaba, Meguro, Tokyo 153-8914, Japan
$^4$ Department of Mathematics, Faculty of Science, Zonguldak B\"ulent Ecevit 
University, Zonguldak 67100, T\"urkiye,\,
e-mail: {\tt myama@ms.u-tokyo.ac.jp}
}



\begin{abstract}
We prove the uniqueness in determining a spatially varying 
zeroth-order coefficient of a one-dimensional time-fractional 
diffusion equation by initial value and Cauchy data at one end point
of the spatial interval.
\end{abstract} 
\baselineskip 18pt

\section{Introduction and main results}
Let $\ell, T>0$.  By $L^2(0,\ell)$, $H^2(0,\ell)$, 
$W^{1,1}(0,T)$, $W^{1,1}(0,T;L^2(0,\ell))$, $L^2(0,T;H^2(0,\ell))$, etc., 
we denote the Lebesgue space and usual Sobolev spaces.
For $0<\alpha<1$, we define the Caputo derivative:
$\ddda v(t) = \frac{1}{\Gamma(1-\alpha)}\int^t_0
(t-s)^{-\alpha} \frac{dv}{ds}(s) ds$ for $v \in W^{1,1}(0,T)$,
where $\Gamma$ denotes the gamma function (e.g., Podlubny \cite{Po}).
Let 
$u = u(x,t) \in W^{1,1}(0,T;\LLLLL) \cap L^2(0,T;H^2(0,\ell))$ satisfy 
$$
\ddda u(x,t) = \ppp_x^2u(x,t) - p(x)u(x,t), \quad
u(x,0) = a(x) \in C^2[0,\ell], \quad 0<x<\ell, \, 0<t<T.
$$ 
This is a model system for example for anomalous diffusion in 
heterogeneous media (e.g., \cite{Po}).
This article is concerned with \\
{\bf Inverse coefficient problem:}
{\it Let $T>0$ be arbitrarily fixed.  
Then determine $p(x)$, $0<x<\ell$ by data
$u(x,0)$ for $0<x<\ell$ and $u(0,t)$, $\ppp_xu(0,t)$
for $0<t<T$.}

We state our first main result on the uniqueness.
\begin{theorem}\label{1}
{\it
Let $p, q\in C^1[0,T], a\in C^2[0,T]$ and
$u, \www{u} \in W^{1,1}(0,T;\LLLLL) \cap L^2(0,T; H^2(0,\ell))$ satisfy
\begin{equation}\label{(1.1)}
\ddda u(x,t) = \ppp_x^2u(x,t) - p(x)u(x,t), \quad u(x,0) = a(x), \quad
\ppp_xu(0, t) = 0, \quad x\in \ell, \, 0<t<T     \end{equation}
and
\begin{equation}\label{(1.2)}
\ddda \www{u}(x,t) = \ppp_x^2\www{u}(x,t) - q(x)\www{u}(x,t), \quad 
\www{u}(x,0) = a(x), \quad
\ppp_xu(0, t) = 0, \quad x\in \ell, \, 0<t<T.     \end{equation}
We assume that 
\begin{equation}\label{(1.3)}
\vert a(x)\vert \ne 0 \quad \mbox{for $0 \le x \le \ell$.}
                              \end{equation}
Then, $u(0,t) = \www{u}(0,t)$ for $0<t<T$, implies 
$p(x) = q(x)$ for $0\le x\le \ell$.
}
\end{theorem}

Theorem \ref{1} is the uniqueness without any data at $x=\ell$, for 
given initial value.  
The condition $\vert a\vert \ne 0$ on $[0,\ell]$ is essential.

The case $\alpha=1$ corresponds to an inverse problem for parabolic 
equations, and we refer to Bukhgeim and Klibanov \cite{BK},
Imanuvilov and Yamamoto \cite{IY98}, Yamamoto \cite{Y09} for example
in general spatial dimensions.
As for the one-dimensional case with $\alpha=1$, see 
Murayama \cite{M}, Pierce \cite{P}, Suzuki \cite{S}, 
Suzuki and Murayama \cite{SM}, where the boundary condition at 
$x=\ell$ is essentially required.  The article Imanuvilov and Yamamoto 
\cite{IY24} proves the uniquness in the case $\alpha=1$ without
data at $x=\ell$.  As for inverse coefficient problems for 
one-dimensional time-fractional diffusion equations with 
$0<\alpha<1$, see Cheng, Nakagawa, Yamamoto and Yamazaki \cite{CNYY} 
as a pioneering work, and 
we can refer for example to 
Jin \cite{J}, Jing and Yamamoto \cite{JY},
Li, Liu and Yamamoto \cite{LLY19}, Liu, Li and Yamamoto
\cite{LiuLY}, Li and Yamamoto \cite{LY} and the references in 
\cite{J}, and we are here limited to 
a few references. 
These articles require the boundary 
condition not only at $x=0$ but also $x=\ell$, which is the same as
\cite{M}, \cite{S}, \cite{SM}.  
In many applications, it is practical to localize data in time only at one end 
point $x=0$, which is our inverse problem.

The key of the proof is the transformation operator (e.g., Levitan 
\cite{L}).  The works \cite{S}, \cite{SM} are based on such an 
operator but the use of the transformation operator in the current article is 
essentially different from them.

In Theorem 1, the assumption of the zero Neumann boundary condition at $x=0$
for solutions $\www u, u$ is restrictive.  We can  drop it 
instead by knowing the values of coefficients $p, q$ near 
$x=0$. 
\begin{theorem}\label{2}
{\it
We assume that $p, q \in C^1[0,\ell], a\in C^2[0,\ell$ and satisfy 
(\ref{(1.3)}). Let $u, \www{u} \in W^{1,1}(0,T;\LLLLL) \cap 
L^2(0,T;H^2(0,\ell))$ satisfy $u(x,0) = \www{u}(x,0) =: a(x)$ for 
$0<x<\ell$ and 
\begin{equation}\label{(1.4)}
\ddda u (x,t)= \ppp_x^2u(x,t) - p(x)u(x,t), \quad 
\ddda \www{u}(x,t) = \ppp_x\www{u}(x,t) - q(x)\www{u}(x,t) \quad 
\mbox{in $(0,\ell) \times (0,T)$}.                   
\end{equation}
We assume that 
there exists $\ep_0\in (0,\ell)$ such that 
$p(x)=q(x)$ for $0\le x\le \ep_0$.
Then, $u(0,t) = \www{u}(0,t)$ and $\ppp_x u(0,t) = \ppp_x \www u(0,t)$ for
$0<t<T$, imply $p(x) = q(x)$ for $0\le x\le \ell$.
}
\end{theorem}


\section{Key lemmata}\label{S2}


\begin{lemma}\label{L1}
{\it
(i)  Let $p, q \in C^1[0,\ell]$. There exists a unique solution $K = K(x,y) \in C^2(\ooo{\OOO})$ 
to the following problem:
\begin{equation}\label{(2.1)}
\left\{ \begin{array}{rl}
& \ppp_x^2{ K}(x,y) -\ppp_y^2 K(x,y)
= q(x) K(x,y) - p(y)K(x,y), \quad (x,y) \in \OOO=\{ (x,y);\, 0<y<x<\ell \}, \\
& (\ppp_y K)(x,0) = 0, \quad  
2\frac{d}{dx}{K}(x,x) = q(x) - p(x),\quad 0<x<\ell, \quad
  K(0,0)=0.
\end{array}\right.
                                   \end{equation}

(ii) {\bf (transformation operator)} 
Let  Let $p, q \in C^1[0,\ell],$ $u \in W^{1,1}(0,T;\LLLLL) \cap L^2(0,T;H^2(0,\ell))$ satisfy 
$\ddda u(x,t) = \ppp_x^2u(x,t) - p(x)u(x,t)$ in $(0,\ell) \times (0,T)$.
Then the function $\www{v}$ given by  
\begin{equation}\label{(2.2)}
\www{v}(t,x) := u(t,x) + \int^x_0 K(x,y)u(t,y) dy, \quad 0<x<\ell, \, 0<t<T
                              \end{equation}
belongs to $W^{1,1}(0,T;\LLLLL) \cap L^2(0,T;H^2(0,\ell))$ and satisfies
\begin{equation}\label{(2.3)}
 \ddda \www{v}(x,t) - \ppp_x^2\www{v}(x,t) + q(x)\www{v}(x,t) 
= -K(x,0)\partial_x u(0,t), \,
 \partial_x\www{v}(0,t) = \partial_x\www{u}(0,t),\quad 
\www{v}(0,t)={u}(0,t).
                                 \end{equation}
}
\end{lemma}
Part (i) is concerned with a Goursat problem and the proof is standard 
by means of the characteristics (e.g., Suzuki \cite{S}).
The part (ii) of Lemma 1 
is the same as Lemma \ref{L2} in \cite{IY24}.

\section{Proofs of Theorems \ref{1} and \ref{2}}\label{S3}

First we show
\begin{lemma}\label{L2} (extension of uniqueness intervals)
{\it
Let functions  $u$, $\www{u} \in W^{1,1}(0,T;\LLLLL) \cap L^2(0,T;H^2(0,\ell))$
satisfy (\ref{(1.4)}).  We assume (\ref{(1.3)}), $p,q\in C^1[0,T]$, $a\in C^2[0,T]$ and 
\begin{equation}\label{(3.1)}
u(0,t) = \www{u}(0,t), \quad \ppp_xu(0,t) = \ppp_x\www{u}(0,t)
\quad \mbox{for $0<t<T$}.               
\end{equation}
(i) We assume that we can find $\delta_0 \in (0,\ell)$  such that 
\begin{equation}\label{(3.2)}
K(x,0)\ppp_xu(0,t) = 0               
\end{equation}
for $0<x<\delta_0$ and $0<t<T$.  Then there exists $x_0 \in (0,\ell)$ 
such that 
\begin{equation}\label{(3.3)}
p(x) = q(x), \quad 0<x<x_0.                            
\end{equation}
(ii) Let (\ref{(3.3)}) holds true with some $x_0\in (0,\ell).$
Moreover we assume that we can find $\delta_1 > x_0$ such 
that (\ref{(3.2)}) holds for $0<x<\delta_1$ and $0<t<T$.
Then, there exists a small constant $\ep > 0$ such that
$p(x) = q(x)$ for $x_0 < x < x_0+\ep$.
}
\end{lemma}
{\bf Proof of Lemma \ref{L2}.}
We define the function $\www{v}$ by (\ref{(2.2)}).
Then, in terms of (\ref{(3.2)}) for $0<x<\delta_0$, 
we can apply Lemma \ref{L1} (ii) to obtain
\begin{equation}\label{(3.5)}
 \ddda \www{v}(x,t) = \ppp_x^2\www{v}(x,t) - q(x)\www{v}(x,t)\quad \mbox{in 
$(0,\delta_0)\times (0,T)$}, 
 \www{v}(0,t) = u(0,t),\quad 
\partial_x\www{v}(0,t) = \ppp_xu(0,t).
\end{equation}
By (\ref{(1.4)}), (\ref{(3.1)}) and (\ref{(3.5)}), the function 
$w:= \www {u} - \www{v}$ satisfies 
$\ddda w(x,t) = \partial^2_x w(x,t) - q(x)w(x,t)$ in $(0,\,\delta_0) \times 
(0,T)$ and $w(0,t) = \ppp_xw(0,t) = 0$ for $0<t<T$.

By the unique continuation for the time-fractional diffusion equation
(Li, Liu and Yamamoto \cite{LLY20}), we see that 
$\www{v}=\www u$ in $(0,\, \delta_0) \times (0,T)$.
Therefore,
\begin{equation}\label{(3.6)}
\int_0^xK(x,y) a(y)dy=0, \quad 0\le x \le \delta_0.        \
\end{equation}
Twice differentiating (\ref{(3.6)}) with respect to $x$, similarly to Lemma 3 
in \cite{IY24}, dividing the corresponding equality  by $a$ and using 
(\ref{(1.3)})
we obtain
\begin{equation}\label{(3.7)}
 (p(x)-q(x))= 
  \int^x_0 (q(x){K}(x,y)-p(y)K(x,y))a(y) dy\,a^{-1}(x)
+ \int^x_0 {K}(x,y)\ppp_y^2 a(y) dy\, a^{-1}(x).
                           \end{equation}
Let us prove a statement (i). Let $\OOO_x := \{ (\xi,\eta);\, 0<\eta<\xi<x\}$.
By the estimate of the solution $K$ to the Goursat problem (\ref{(2.1)}) 
for any $x\in [0,\ell]$,
we have
\begin{equation}\label{(3.9)}
\Vert K\Vert_{C(\ooo{\OOO_x})} \le C \Vert p-q\Vert_{C[0,x]}.
                                  \end{equation}
 From (\ref{(3.7)}), (\ref{(3.9)}) we have 
$$
\Vert p-q\Vert_{C[0,x]}\le C x\Vert p-q\Vert_{C[0,x]}.
$$
Taking $x$ sufficiently small, we obtain (\ref{(3.3)}).

Next we prove (ii) provided that (\ref{(3.3)}) holds with 
$x_0 \in (0,\delta_0)$.
Since (\ref{(3.3)}) implies $K(x,x) = 0$ for $0\le x\le x_0$, the uniqueness of
solution to (\ref{(2.1)}) yields
\begin{equation}\label{(3.4)}
K(x,y) = 0 \quad \mbox{if $0\le y \le x\le x_0$}.    
\end{equation}
Then, we can follow the proof of Lemma 3 in \cite{IY24} and we provide the 
proof for completeness.

We divide the two integral terms on the right-hand side of (\ref{(3.7)}):
\begin{align*}
& \int^x_0 (q(x){K}(x,y)-p(y)K(x,y))a(y) dy\,a^{-1}(x)
+ \int^x_0 {K}(x,y)\ppp_y^2 a(y) dy\, a^{-1}(x)=\\
& \left( \int^x_{x_0} + \int^{x_0}_0\right)
\{ (q(x){K}(x,y)-p(y)K(x,y))a(y) + K(x,y)\ppp_y^2 a(y)\} dy\, a^{-1}(x)
=: I_1(x) + I_2(x), \quad x_0 \le x \le \delta_0.
\end{align*}
Then we have
\begin{equation}\label{(3.8)}
\vert I_1(x)\vert \le C\vert x-x_0\vert \sup_{(\xi,\eta)\in \OOO_x}
\vert K(\xi,\eta)\vert \quad \mbox{for $x_0\le x\le \delta_0$.}                                                   
\end{equation}
Here and henceforth $C>0$ denotes generic constants depending on 
$a, p, q$ but not on $x, x_0$.

By (\ref{(3.8)}) and (\ref{(3.9)}), we obtain
\begin{equation}\label{(3.10)}
\vert I_1(x)\vert \le C\vert x-x_0\vert \Vert p-q\Vert_{C[0,x]}
\quad \mbox{for $x_0\le x\le \delta_0$}.
                                         \end{equation}
Moreover by (\ref{(3.3)}), we can represent as follows:
\begin{align*}
& I_2(x)=
 \int^{x_0}_0 ({K}(x,y)-K(x_0,y))\ppp_y^2 a(y) dy\, a^{-1}(x)
\end{align*}
for all $x\in (x_0, \delta_0)$.
Then, the mean value theorem yields
\begin{equation}\label{(3.11)}
\vert I_2(x)\vert \le Cx_0 \sup_{y\in [0,x_0]}\vert K(x,y)-K(x_0,y)\vert
\le C\sup_{(z,y)\in [0,x]\times [0,x_0]}\vert \partial_z K(z,y)\vert 
\vert x-x_0\vert.               
\end{equation}
From (\ref{(3.7)}), (\ref{(3.10)}) and (\ref{(3.11)}), we obtain
$\vert (p-q)(x)\vert \le C\vert x-x_0\vert \Vert p-q\Vert_{C[0,x]}$,
that is, $\Vert p-q\Vert_{C[0,x]}
\le C \vert x-x_0\vert \Vert p-q\Vert_{C[0,x]}$ for $x_0 \le x \le \delta_0$.
Taking sufficiently small $\ep>0$, we obtain 
$\Vert p-q\Vert_{C[0,x_0+\ep]}=0$.  


Thus the proof of Lemma \ref{L2} is complete. 
$\blacksquare$

{\bf Proof of Theorem \ref{1}.}
Let $x_0 \in (0,\ell)$ be the maximum value satisfying 
$p(x) = q(x)$ for $0\le x\le x_0$.  
By the assumption $\ppp_xu(0,t) = 0$ for $0<t<T$,
condition (\ref{(3.2)}) is satisfied for $0<x<\ell$.  
In view of Lemma \ref{L2} (i), such an 
$x_0>0$ exists.
If $x_0 = \ell$, then the proof is already finished.
Therefore, we can assume that $0<x_0<\ell$. 
Then, the statement (ii) of the Lemma \ref{L2} implies the 
existence of $\ep>0$ such that $p=q$ in $(0, x_0+\ep)$.  This contradicts
the definition of $x_0$.  Thus the proof of Theorem \ref{1} is complete.
$\blacksquare$

{\bf Proof of Theorem \ref{2}.} Let $x_0 \in (0,\ell)$ be the maximum value 
satisfying 
$p(x) = q(x)$ for $0\le x\le x_0$.  Theorem \ref{2} implies $x_0\ge\ep_0>0$.
Then equality (\ref{(3.4)}) holds true. 
The uniqueness of the  Cauchy problem for the  hyperbolic equation implies 
that $K(x,0)=0$ on  $\delta_0=\min \{2\ep_0,\ell\}$.
Then (\ref{(3.2)}) holds true.
Lemma \ref{L1} (ii) yields $p(x) = q(x)$ for $0\le x\le x_0+\ep_0$. 
Since $\ep_0>0$, we reach a contradiction. 
$\blacksquare$

\section{Concluding remarks}\label{S4}

We proved the uniqueness for the inverse coefficient 
problems without full boundary conditions 
for spatially one-dimensional time-fractional diffusion equations.
Our method is based on the transformation operator, and 
widely applicable if we have the unique continuation property 
for the partial differential equations under consideration.
Our approach requires that the spatial dimension is one.

{\bf Acknowledgments.}
The work was supported by Grant-in-Aid for Scientific Research (A) 20H00117 
and Grant-in-Aid for Challenging Research (Pioneering) 21K18142 of 
Japan Society for the Promotion of Science.

\end{document}